\documentclass[12pt]{article}
\usepackage{latexsym}
\usepackage{amssymb}
\usepackage{amsthm}
\begin{document}
\newtheorem{theorem}{Theorem}[section]
\newtheorem{lemma}[theorem]{Lemma}
\newtheorem{proposition}[theorem]{Proposition}
\newtheorem{corollary}[theorem]{Corollary}
\newtheorem*{definition}{Definition}
\newtheorem*{question}{Question}
\newtheorem{conjecture}{Conjecture}
\newtheorem*{thm}{Theorem}
\newcommand{\F}{\ensuremath{\mathbb F}}
\newcommand{\N}{\mathcal N}
\newcommand{\R}{\mathcal R}
\newcommand{\Z}{\mathbb Z}

\title{Conjugacy classes of maximal cyclic subgroups}
\author{M. Bianchi, R.D. Camina, Mark L. Lewis \& E. Pacifici}
\date \
\maketitle

\begin{abstract} In this paper, we set $\eta (G)$ to be the number of conjugacy classes of maximal cyclic subgroups of $G$.  We consider $\eta$ and direct and semi-direct products. We characterize the normal subgroups $N$ so that $\eta (G/N) = \eta (G)$.  We set $G^- = \{ g \in G \mid \langle g \rangle {\rm ~is~not ~maximal~cyclic} \}$.  We show if $\langle G^- \rangle < G$, then $G/\langle G^- \rangle$ is either (1) an elementary abelian $p$-group for some prime $p$, (2) a Frobenius group whose Frobenius kernel is a $p$-group of exponent $p$ and a Frobenius complement has order $q$ for distinct primes $p$ and $q$, or (3) isomorphic to $A_5$.  \\[1ex]
{\it Keywords:} direct product, exponent $p$, covering\\[1ex]
{\it 2020 Mathematics Subject Classification:} 20D15, 20E34
\end{abstract}

\section{Introduction}

Unless otherwise stated, all groups in this paper are finite.  A {\it covering} of a group $G$ is a set of
proper subgroups $\{ H_i\}$, called {\it components}, such that $G \subseteq \bigcup_i H_i$.  A covering is called {\it irredundant} if removing any component means the set is no longer a covering.  Identifying coverings of groups and the minimal size of such a covering, called the covering number, has a long mathematical history.  It is an easy problem to show that a group cannot be covered by $2$ subgroups.  In 1926, Scorza considered groups with a covering of size $3$, \cite{scorza}; then, in 1994, Cohn continued the investigation, considering groups with a minimal covering of sizes $4$, $5$, and $6$ amongst other things, \cite{cohn}.  The case of $7$ is considerably more difficult and surprising in that no such covering exists.  This was conjectured by Cohn and proved by Tomkinson.  Tomkinson also determined the coverings of solvable groups in \cite{tomkinson}.   Since then, many authors have written on the theme; for an overview see \cite{bhargava} or \cite{Kappe_survey}.  At this time, the research appears to be on determining the coverings of nonabelian simple groups and the symmetric groups.

Recently, authors have considered {\it normal coverings}, that is coverings which are invariant under $G$-conjugation (not coverings by normal subgroups). The normal covering number is the smallest number of conjugacy classes of proper subgroups in a normal covering of $G$.   In \cite{praeger}, \cite{praeger2} and \cite{praeger3}, the authors consider the normal covering number for  symmetric and alternating groups, while in \cite{garonzi} the authors consider groups for which the normal covering number is equal to the covering number.

Also, there has been interest in coverings where certain requirements are imposed on the components of the covering.  In this paper, we consider coverings by cyclic subgroups up to conjugacy.  For finite groups, it is not difficult to see that every cyclic subgroup is contained in a maximal cyclic subgroup.  Furthermore, every covering by cyclic subgroups has to contain the maximal cyclic subgroups since the only cyclic subgroups containing the generators of the maximal cyclic subgroups will be the maximal cyclic subgroups.  Therefore, the set of the maximal cyclic subgroups of $G$ form the only irredundant covering of $G$ by cyclic subgroups.   Hence, it makes sense to study the set of maximal cyclic subgroups of $G$.   In particular, the set of conjugacy classes of maximal cyclic subgroups form an irredudant normal covering where all components are cyclic.  

We note that in infinite groups cyclic subgroups need not be contained in maximal cyclic subgroups (see the Pr\"ufer group) so the question of irredundant coverings of infinite groups by cyclic subgroups may be independent of the question of maximal cyclic subgroups.  This is being explored by the second and third authors with Yiftach Barnea and Mikhail Ershov in \cite{pre}.

Coverings of abelian groups by maximal cyclic subgroups have been studied by Rog\'erio in \cite{rogerio}.  For nonabelian groups, von Puttkamer initiated the study of the number of conjugacy classes of maximal cyclic subgroups in Section 5 of his dissertation (\cite{von}).  In particular, von Puttkamer considers finite groups $G$ that have two conjugacy classes of maximal cyclic subgroups.  He proves that such groups must be solvable, and in addition, he proves that they have derived length at most $4$.  In fact, he essentially classifies them.  We pick up the study of groups in terms of the number of conjugacy classes of maximal cyclic subgroups in this paper.   With this in mind, we set the following definition.   

\begin{definition} Let $G$ be  finite group, we denote the number of conjugacy classes of maximal cyclic subgroups of $G$ by $\eta(G)$.
\end{definition}



In Section \ref{prod sect}, we will consider $\eta (G)$ when $G$ is a direct product.  We also determine $\eta (G)$ when $G$ is a Frobenius group.  For a group $G$, we define $G^-$ to be the set of $g \in G$ such that $\langle g \rangle$ is not maximal cyclic.  In Section \ref{eta and quo}, we will show that there is a strong connection between $\eta (G)$ and $G^-$.  In Section \ref{G^- as set}, we consider some properties of $G^-$ as a set.

In particular, we prove the following:

\begin{theorem}\label{1st main}
If $G$ is a group and $\langle G^- \rangle < G$, then $G/\langle G^- \rangle$ is either (1) a $p$-group with exponent $p$, (2) a Frobenius group whose Frobenius kernel is a $p$-group with exponent $p$ and whose Frobenius complements have order $q$ for some prime $q \ne p$, or (3) $A_5$.
\end{theorem}

We will see that we have stronger control on $G^-$ when $G$ is a $p$-group.  When $G$ is a $p$-group and $N$ is a normal subgroup so that $G/N$ has exponent $p$, then $G^- \subseteq N$.

Let $N$ be a normal subgroup of $G$ then it is not hard to see that $\eta(G/N) \leq \eta(G)$. In Theorem \ref{quot} we give criteria for when equality occurs.  We note that the criteria is in terms of $G^-$.  We then show that when $G$ is a $p$-group, there exists a characteristic subgroup $X(G)$ of $G$ such that $\eta(G) = \eta(G/X(G))$ and $X(G)$ is maximal under this condition.  We will give an example to see that we cannot find such a subgroup in general.

We close by proving a couple of results regarding the structure of $N$ when $\eta (G/N) = \eta (G)$.  The first considers the relationship between $N$ and the derived subgroup, $G'$.

\begin{theorem}\label{3rd main}
Suppose $G$ is a group and $N$ is a subgroup so that $\eta (G) = \eta (G/N)$.  If every nontrivial Sylow subgroup of $G/G'$ is noncyclic, then $N \le G'$.
\end{theorem}

Finally, when $G$ is a $p$-group and a maximal cyclic subgroup is normal in $G$, then we can prove that $X(G)$ is cyclic.

\begin{theorem}\label{4th main}
Suppose $G$ is a $p$-group for some prime $p$.  If some maximal cyclic subgroup of $G$ is normal, then $X(G)$ is cyclic.
\end{theorem}

\section{$\eta$ and products} \label{prod sect}
 
In this section, we consider the relationship between $\eta$ and direct products and semi-direct products.  In this next result we give some idea of the relationship between $\eta$ and direct products.  
 
 
 
\begin{lemma} \label{dir product}
Suppose $G = H \times K$, then the following are true:
\begin{enumerate}
\item [(i)] $\eta(G) \ge \eta(H)\eta(K)$.
\item [(ii)] If $(|H|,|K|) = 1$, then $\eta(G) = \eta(H)\eta(K)$.
\item [(iii)] If $p$ is a prime so that $p$ divides $|K|$ and $H$ is a nontrivial $p$-group, then $\eta (G) \ge \eta (H) \eta (K) + \eta_p (K) > \eta (H) \eta (K)$ where $\eta_p (K)$ is the number of conjugacy classes of maximal cyclic subgroups of $K$ whose order is divisible by $p$.
\item [(iv)] If $H$ is nilpotent and $p$ divides both $|H|$ and $|K|$, then $\eta (G) > \eta (K)$.
\item [(v)] If $H$ and $K$ are both nontrivial $p$-groups, then $\eta (G) \ge \eta (H)\eta (K) + \eta (H) + \eta (K)$. 
\end{enumerate} 
\end{lemma}
 
{\bf Proof.} Suppose $\langle h \rangle$ is a maximal cyclic subgroup of $H$ and $\langle k \rangle$ is a maximal cyclic subgroup of $K$.  Suppose that $\langle (h,k) \rangle \le \langle g \rangle$ for some element $g \in G$.  We can write $g = (h_1,k_1)$ for $h_1 \in H$ and $k_1 \in K$.  Since $(h,k) \in \langle g \rangle$, there is an integer $n$ so that $g^n = (h,k)$, and hence, $h_1^n = h$ and $k_1^n = k$.  It follows that $\langle h \rangle \le \langle h_1 \rangle$ and $\langle k \rangle \le \langle k_1 \rangle$.  By maximality, we have $\langle h \rangle = \langle h_1 \rangle$ and $\langle k \rangle = \langle k_1 \rangle$, and so $\langle g \rangle = \langle (h,k) \rangle$.  Thus, $\eta (H) \eta (K) \le \eta (G)$.   
 
We now assume that $(|H|,|K|) = 1$.  Let $\langle g \rangle$ be a maximal cyclic subgroup of $G$.  We know $g = (h,k)$ for unique $h \in H$ and $k \in K$.  Suppose $\langle h \rangle \le \langle h_1 \rangle$ for some $h_1 \in H$ and $\langle k \rangle \le \langle k_1 \rangle$ for some $k_1 \in K$.  Since $o(h_1)$ and $o(k_1)$ are coprime, it follows that $(h,1)$ and $(1,k)$ are powers of $(h_1,k_1)$.  Thus, $g = (h,1) (1,k) \in \langle (h_1,k_1) \rangle$, and so, $\langle g \rangle \le \langle (h_1,k_1) \rangle$.  Maximality of $\langle g \rangle$ implies that $\langle g \rangle = \langle (h_1,k_1) \rangle$, and so, $\langle h \rangle = \langle h_1 \rangle$ and $\langle k \rangle = \langle k_1 \rangle$.  We conclude that $\langle h \rangle$ and $\langle k \rangle$ are maximal cyclic subgroups.  This implies that $\eta (G) \le \eta (H) \eta (K)$, and so we have $\eta (G) = \eta (H) \eta (K)$.   
 
Next, we have that $H$ is a nontrivial $p$-group and $p$ divides $|K|$.  Consider an element $k \in K$ so that $\langle k \rangle$ is maximal cyclic in $K$ and $p$ divides $o (k)$.  We claim that $\langle (1,k) \rangle$ will be a maximal cyclic subgroup of $H \times K$.  Suppose there exists an element $g \in G$ so that $g^q = (1,k)$ for some prime $q$ that divides $o(g)$.  We can write $g = (h^*,k^*)$ for $h^* \in H$ and $k^* \in K$.  We see that $(1,k) = g^q = ((h^*)^q,(k^*)^q)$ and so, $(k^*)^q = k$.  We know that $o(g)$ is the least common multiple of $o(h^*)$ and $o (k^*)$.  Notice that $o(k)$ divides $o(k^*)$ and so, $p$ divides $o(k^*)$.  Since $o(h^*)$ is a power of $p$, we see that $q$ divides $o(k^*)$ and this contradicts the choice of $\langle k \rangle$ as a maximal cyclic subgroup of $K$.  Notice that $(1,k)$ is not conjugate to any of the subgroups included in the count in part (i).  This yields $\eta (G) \ge \eta (H) \eta (K) + \eta_p (K)$.  By Cauchy's theorem, we know that $K$ has an element whose order is divisible by $p$, and so, $K$ has a maximal cyclic subgroup whose order is divisible by $p$, and so, $\eta_p (K) \ge 1$.  This proves (iii).
 
We now prove conclusion (iv).  We can write $H = P \times Q$ where $P$ is the Sylow $p$-subgroup of $H$ and $Q$ is the Hall $p$-complement.  We then have $G = (K \times Q) \times P$.  By conclusion (iii), we have $\eta (G) > \eta (K \times Q)$ and by conclusion (i), we have $\eta (K \times Q) \ge \eta (K)$.  Combining these, we have conclusion (iv).  For conclusion (v), notice that every maximal cyclic subgroup in either $H$ or $K$ will have order divisible by $p$.  Thus, if $h \in H$ and $k \in K$ satisfy that $\langle h \rangle$ is a maximal cyclic subgroup of $H$ and $\langle k \rangle$ is a maximal cyclic subgroup of $K$, then arguing as in (iii), we see that $\langle (h,1) \rangle$ and $\langle (1,k) \rangle$ are maximal cyclic subgroups of $G$.  Since these cannot be conjugate to the maximal cyclic subgroups counted in (i), we obtain $\eta (G) \ge \eta (H) \eta (K) + \eta (H) + \eta (K)$.   
$\Box$\\[1ex]
 
We note that equality can hold in Lemma \ref{dir product} even if $|H|$ and $|K|$ are not relatively prime.  Take $H = S_3$ and $K$ to be the dihedral group of order $10$.  Then $\eta (H \times K) = 4$ and $\eta (H) = \eta (K) = 2$.  The various cases that arise show that computing $\eta$ for direct products is not trivial when the factors do not have coprime orders.
 
We expect the situation for computing $\eta$ of semi-direct products to be even more complicated.  We begin with the case of Frobenius groups which is relatively straightforward.
 
\begin{proposition}
Let $G = NH$ be a Frobenius group with Frobenius kernel $N$ and Frobenius complement $H$.  Then $\eta (G) = \eta^* (N) + \eta (H)$ where $\eta^* (N)$ is the number of $H$-orbits on the $N$-conjugacy classes of maximal cyclic subgroups of $N$.
\end{proposition}
 
{\bf Proof:}
Let $C$ be a maximal cyclic subgroup of $G$.  We know that either $C \le N$ or some conjugate of $C$ lies in $H$.  Furthermore, if $C \le N$, then $C$ is a maximal cyclic subgroup of $G$ and if $C \le H$, then $C$ is maximal cyclic in $H$.  Also, we know that if $C$ and $C^g$ lie in $H$, then $1 < H \cap H^g$ and so, $g \in H$.  It follows that the number conjugacy classes of maximal cyclic subgroups of $G$ will equal the number of $H$-conjugacy classes of maximal cyclic subgroups in $H$ plus the number of $H$-orbits on the $N$-conjugacy classes of maximal cyclic subgroups of $N$.
$\Box$\\
 
In general, when $G$ is a semi-direct product of $H$ acting on $N$ where $(|N|,|H|) = 1$, we would like to bound $\eta (G)$ in terms of $\eta^* (N)$ and $\eta (H)$.  It is tempting to ask whether $\eta^* (N) + \eta (H)$ is the lower bound for $\eta (G)$.  It turns out the answer is no.  Let $G$ be the semidirect product of $D = D_8$ acting faithfully on $N = Z_3 \times Z_3$.   (In the SmallGroups library in Magma, this is SmallGroup (72,50).)  It is not difficult to see that $D$ has two orbits on cyclic subgroups of order $3$ in $N$ and so, $\eta^* (N) = 2$.  We know $\eta (D) = 3$; so $\eta^* (N) + \eta (D) = 5$.  On the other hand, one can see that the maximal cylic subgroups of $G$ have orders $4$ and $6$, and that there is one conjugacy class of cyclic subgroups of order $4$ and two conjugacy classes of cyclic subgroups of order $6$; so $\eta (G) = 3$.  We now provide a lower bound; however, we do not have any examples where this lower bound is met. 
 
\begin{proposition}\label{centre}  
Let $N$ be a normal subgroup of a group $G$ and let $\eta^*(N)$ be the number of $G$-orbits on the $N$-conjugacy classes of maximal cyclic subgroups of $N$.  Then $\eta (G) \geq \eta^*(N)$.  In particular,
\begin{enumerate}
\item [(i)] if $N$ is central in $G$ then $\eta(G) \geq \eta(N)$.
\item [(ii)] if $|G:N| = k$ then $\eta (G) \geq \eta(N)/k$.
\end{enumerate}
\end{proposition}
 
{\bf Proof.}   
Let  $C_1, \ldots, C_m$ be a set of representatives for the $G$-orbits of the maximal cyclic subgroups of $N$.  We know that each $C_i$ is contained in some maximal cyclic subgroup of $G$.  It suffices to show that if $i \ne j$, then $C_i$ and $C_j$ are not $G$-conjugate to subgroups of the same maximal cyclic subgroup of $G$.  Without loss of generality, we may assume $i = 1$ and $j = 2$.  Suppose that $C_1$ and $C_2$ are $G$-conjugate to subgroups of the same maximal cyclic subgroup $D$ of $G$. Replacing $C_1$ and $C_2$ by conjugates if necessary, we may assume that $C_1$ and $C_2$ are subgroups of $D$.  Since $C_1$ and $C_2$ are contained in $N$, we see that $C_1$ and $C_2$ are subgroups of $D \cap N$.  Now, $D \cap N$ is a cyclic subgroup of $N$ and $C_1$ and $C_2$ are maximal cyclic subgroups of $N$; so we have $C_1 = D \cap N = C_2$, but this contradicts the fact that $C_1$ and $C_2$ are representatives of distinct orbits.
$\Box$\\

We close by mentioning that the situation for subgroups is more complicated. Consider the wreath product $C_3 \wr C_3$ and its base group $C_3 \times C_3 \times C_3$. Then $\eta(C_3 \wr C_3) = 7$ but $\eta(C_3 \times C_3 \times C_3) = 13$.  

\section{The set $G^-$}\label{G^- as set}

Let $G$ be a finite group. We define the set $$G^- = \{ g \in G: \langle g \rangle \; {\rm is \; not\; a \; maximal \; cyclic \; group}\}.$$  For a prime $p$, we set $G^{\{p\}} = \{g^p: g \in G \}$.

We now give a characterization of this set.

\begin{lemma} \label{pgrp}
Let $G$ be a group.  Then 
$$G^- = \{ g^q \mid g \in G,  q {\rm \; a\; prime\; dividing\; the \; order\; of \;} g \}.$$
In particular, if $G$ is a $p$-group then $G^- = G^{\{p\}}$.
\end{lemma}

{\bf Proof:}
Suppose $x \in \{ g^q \mid g \in G, q {\rm ~a ~prime~and~} q | o(g) \}$.  There exists $y \in G$ so that $y^q = x$ for some prime $q$ with $q$ dividing $o (y)$.  We see that $\langle x \rangle < \langle y \rangle$, and thus, $x \in G^-$.  This yields $\{ g^q \mid g \in G, q {\rm ~a ~prime~and~} q | o(g) \} \subseteq G^-$.  Now, suppose $x$ is in $G^{-}$, so that $\langle x \rangle$ is a proper subgroup of $\langle h \rangle$ for some $h \in G$, so some prime $q$ must divide $|\langle h \rangle:\langle x \rangle |$. We can then choose an element $y$ of $\langle h \rangle$ such that $\langle y \rangle$ contains $\langle x \rangle$ and $|\langle y \rangle : \langle x \rangle | = q$, so $\langle x \rangle = \langle y^q \rangle$ with $q$ dividing $o(y)$.  In particular, $o (x) = o (y)/q$.  We know that $x = (y^q)^l$ for some integer $l$ that is relatively prime to $o(y^q) = o (y)/q$.  If $q$ does not divide $l$, then $l$ is relatively prime to $o(y)$.  We have that $x = y^{ql} = (y^l)^q$ and $q$ divides $o(y^l) = o (y)$.  We now suppose that $q$ does divide $l$.  This implies that $q$ does not divide $m = o(y)/q = o(y^q)$.  Set $l^* = l + m$.  Hence, $(y^q)^{l^*} = (y^q)^{l + m} = y^{ql} = x$.  We prove that $l^*$ is coprime to $o(y)$.  Since $q$ divides $l$ and not $m$, we see that $q$ does not divide $l^*$.  Suppose $p$ is a prime divisor of $m = o(y)/q$.  Thus, $p$ does not divide $l$; so again $p$ does not divide $l^*$.  It follows that $l^*$ is coprime to $o(y)$.  We have $x = (y^{l^*})^q$ and $q$ divides $o(y^{l^*})$.  This proves the claim.  In particular, we have $G^- \leq \{ g^q \mid g \in G, q {\rm ~a ~prime~and~} q | o(g) \}$, and equality follows.
$\Box$\\

In this next result, we show that when $G$ is a $p$-group, $G^-$ is contained in every normal subgroup whose quotient has exponent $p$.  We also show that it is rare when $G$ is not a $p$-group for $G^-$ to be contained in a proper subgroup.  Notice that Theorem \ref{1st main} follows from conclusion (1).  

\begin{theorem}
Let $G$ be a group.  Then the following statements are true:
\begin{enumerate}
\item If $N$ is a normal subgroup of $G$, then $G^- \subseteq N$ if and only if every element in $G \setminus N$ has prime power order and every element in $G/N$ has prime order.  In particular, $G/N$ is either (1) a $p$-group with exponent $p$, (2) a Frobenius group whose Frobenius kernel is a $p$-group with exponent $p$ and whose Frobenius complement has order $q$ for some prime $q \ne p$, or (3) $A_5$.
\item If $M$ is a normal subgroup of $G$ and $p$ is a prime so that $|G:M| = p$, then $G^- \subseteq M$ if and only if every element in $G \setminus M$ has $p$-power order.
\item If $G$ is a $p$-group and $N$ is a normal subgroup so that $G/N$ has exponent $p$, then $G^- \subseteq N$.
\end{enumerate}
\end{theorem}

{\bf Proof:}
We begin by considering an element $g \in G$ whose order is divisible by at least two primes, say $p$ and $q$.  It follows that $g^p$ and $g^q$ lie in $G^-$.  It is not difficult to see that $\langle g \rangle = \langle g^p, g^q \rangle$, so $g \in \langle G^- \rangle$.  It follows that if $\langle G^- \rangle \le N$ where $N$ is a normal  subgroup of $G$, then every element in $G \setminus N$ must have prime power order.  Furthermore, if $p$ is the prime that divides $o(g)$, then $g^p \in G^- \subseteq N$ and so, $(gN)^p = N$.  Hence, every element of $G/N$ has prime order.   Conversely, suppose $N$ is a normal subgroup of $G$ so that every element in $G \setminus N$ has prime power order and every element in $G/N$ has prime order.  Let $h \in G^-$.  We know that $h = g^p$ where $g \in G$ and $p$ is a prime that divides $o (g)$.  If $o (g)$ is also divisible by some prime other than $p$, then we cannot have $g \in G \setminus N$, so $g$ and hence $h$ must lie in $N$.  We may now suppose that $o(g)$ is a power of $p$.  If $g \in N$, then we are done, so we assume $g \not\in N$.  We see that $o(gN)$ divides $o(g)$ and is prime, we have $o(gN) = p$ and so, $h = g^p \in N$.  This proves that $G^- \subseteq N$. 

The study of groups where all the elements have prime order has a long history.  It appears that they were first studied by Higman in \cite{Higman2}.  This work continued in \cite{brandl}, \cite{Deaconescu89}, \cite{ShiYang}, \cite{ShiYang2}, and \cite{suzuki}.  The study of such groups culminated by being completely classified by Cheng, et. al in \cite{Cheng93}. In particular, they prove that if $G$ is a group where all elements have prime order, then $G$ is one of the following: (i) a $p$-group of exponent $p$, (ii) a Frobenius group whose Frobenius kernel is a $p$-group of exponent $p$ and whose Frobenius complement has order $q$ for distinct primes $p$ and $q$, or (iii) $A_5$.  Applying this result to $G/N$ completes the proof of (1).

We have the condition that $M$ is normal in $G$ and $|G:M| = p$.  Note that if every element in $G \setminus M$ has $p$-power order, then (1) implies that $G^- \subseteq M$.  On the other hand, if $G^- \subseteq M$, then every element $g \in G \setminus M$ has prime power order by (1) and $o (gM) = p$ which implies that $p$ divides $o(g)$.  It follows that every element in $G \setminus M$ has $p$-power order.  Suppose $G$ is a $p$-group and $N$ is a normal subgroup so that $G/N$ has exponent $p$.  Then we may apply (1) to see that $G^- \subseteq N$.  
$\Box$\\

We note that there are $p$-groups $G$ where $G^{\{p\}}$ is a subgroup of $G$ and $p$-groups where $G^{\{p\}}$ is not a subgroup of $G$.  We next show that for $G^-$ to be a subgroup, it must be the case that all of the elements of $G$ have prime power order.

\begin{lemma}
Let $G$ be a group.  If $G^-$ is a subgroup of $G$, then every element of $G$ must have prime power order.
\end{lemma}

{\bf Proof:}
Suppose $h$ is an element of $G$ that does not have prime power order.  Then there is an element $g$ so that $\langle g \rangle$ is maximal cyclic and $h \in \langle g \rangle$.  It follows that $o(g)$ is not a prime power.   Hence, there exist nonidentity elements $a, b$ that are powers of $g$ so that $g = ab$ and $(o(a), o (b)) = 1$.  We see that $g \not\in G^-$, but $a,b \in G^-$, and this implies that all elements of $G$ have prime power order.
$\Box$ \\

Higman in \cite{Higman2} classifies solvable groups where all elements have prime power order.  However, we present a different viewpoint to describe these groups.  In particular, if $G$ is a solvable group, then all the elements have prime power order if and only if (1) $G$ is a $p$-group or (2) $G$ is either a Frobenius group or a $2$-Frobenius group that is also a $\{p, q\}$-group for primes $p$ and $q$.  Recall that a group $G$ is a $2$-Frobenius group if there exist normal subgroups $L \le K \le G$ so that $G/L$ and $K$ are Frobenius groups with Frobenius kernels $K/L$ and $L$, respectively.  This can be seen using the {\it Grueneberg-Kegel graph (GK-graph)} (sometimes called the prime graph).  This graph has the primes dividing $|G|$ as its vertices, and there is an edge between $p$ and $q$ if there is some element $g \in G$ so that $pq$ divides $o (g)$.   The GK-graph of a group where all the elements have prime power order will be an empty graph (i.e., a graph with no edges).  Note that the GK-graph of $G$ consists of a single vertex if and only if $G$ is a $p$-group.  When $G$ is solvable, it is known that the GK-graph is disconnected if and only if $G$ is a Frobenius group or a $2$-Frobenius group and that it has two connected components (see \cite{williams}).  Assuming all elements of $G$ have prime power order, this implies that $|G|$ is divisible by two primes $p$ and $q$.

Let $p$ and $q$ be primes so that $p$ divides $q-1$.  It is not difficult to see that there will exist a Frobenius group $G$ whose Frobenius kernel is cyclic of order $q^2$ and a Frobenius complement has order $p$.  It is easy to see that $G^-$ is going to be the unique subgroup of order $q$.  This gives an example of a group $G$ where $G^-$ is a subgroup.

Now let $p$ and $q$ be primes so that $p^2$ divides $q-1$.  In this case, we take $G$ to be the Frobenius group whose Frobenius kernel has order $q$ and whose Frobenius complement has order $p^2$.  It is not difficult to see in this case that $G^-$ will be the elements of order $p$ in $G$ along with $1$, and that they do not form a subgroup.  This gives an example where $G^-$ is not a subgroup.

Let $l(G)$ be the number of conjugacy classes of cyclic subgroups of $G$.  Notice that $\eta (G) \leq l (G) - 1$ when $G$ is a nontrivial group.  If every element of $G$ has prime order, then every cyclic subgroup will be maximal cyclic and so $\eta (G) = l (G) - 1$.  Conversely, if $l(G) - 1 = \eta (G)$, then every nontrivial cyclic subgroup will be maximal cyclic and so will have prime order.  It will follow that every element of $G$ will have prime order.  Hence, $\eta (G) = l (G) - 1$ if and only every element of $G$ has prime order.  Note that we have seen that every element of $G$ has prime order if and only if (1) $G$ is a $p$-group of exponent $p$, (2) $G$ is a Frobenius group where the Frobenius kernel is a $p$-group of exponent $p$ and a Frobenius complement has order $q$ for primes $p$ and $q$, or (3)  $G \cong A_5$.

\section{$\eta$ and quotients} \label{eta and quo}


In this next result we show that $\eta (G/N) \leq \eta (G)$ when $N$ is a normal subgroup.  We then give necessary and sufficient conditions for $\eta (G/N) = \eta (G)$.  We also refine the necessary and sufficient conditions when $G$ is a $p$-group.  Let $S$ and $T$ be subsets of $G$, then we set $ST = \{ st \mid s \in S, t \in T\}$.  When $N$ is a normal subgroup of $G$, we write $SN/N = \{ sN \mid s \in S \}$ as a subset of $G/N$.

\begin{theorem}\label{quot} 
Let $N$ be a proper normal subgroup of $G$. Then 
\begin{enumerate}
\item $\eta(G/N) \leq \eta(G)$.
\item $\eta(G/N) = \eta(G)$ if and only if the following conditions hold:
\begin{enumerate}
\item $N\subseteq G^-$,
\item $(G/N)^- = \{gN \in G/N \mid gN \subseteq G^-\}$, and
\item Every element $x \in G \setminus G^-$ satisfies the condition that every element of $xN \setminus G^-$ is conjugate to a generator of $\langle x \rangle$.
\end{enumerate}
\item  If $G$ is a $p$-group and $\eta(G/N) = \eta(G)$, then $G^-$ is a union of $N$-cosets.
\item $\eta (G/N) = \eta(G)$ and $G^-$ is a union of $N$-cosets if and only if for every $x \in G \setminus G^-$ every element of $xN$ is conjugate to a generator of $\langle x \rangle$.
\item  If $\eta (G/N) = \eta(G)$ and $G^-$ is a union of $N$-cosets, then the following hold: 
\begin{enumerate}
\item $G^-N = G^-$, and 
\item $(G/N)^- = G^-N/N$.
\end{enumerate}   
\end{enumerate}
\end{theorem}

{\bf Proof.} 
Fix an element $x \in G$.  First note $\langle x \rangle$ is maximal cyclic if and only if $x \in G \setminus G^-$.  Similarly, $\langle xN \rangle$ is maximal cyclic if and only if $xN \in G/N \setminus (G/N)^-$.  
 

Suppose $\langle xN \rangle$ is maximal cyclic in $G/N$.  Let $\langle g \rangle$ be a maximal cyclic subgroup containing $\langle x \rangle$.  It follows that $\langle xN \rangle \le \langle gN \rangle$.  By the maximality of $\langle xN \rangle$, we have that $\langle xN \rangle  = \langle gN \rangle$.  Let $x_1, \dots, x_n \in G$ be chosen so that $\langle x_1 N \rangle, \dots, \langle x_n N \rangle$ are representatives of the conjugacy classes of maximal cyclic subgroups of $G/N$.  For each $i$, choose $g_i \in G$ so that $\langle g_i \rangle$ is a maximal cyclic subgroup of $G$ containing $\langle x_i \rangle$.   The work above shows that $\langle x_i N \rangle = \langle g_i N \rangle$.  Suppose $\langle g_i \rangle$ and $\langle g_j \rangle$ are conjugate for $1 \le i, ~j \le n$.  Hence, there exists $y \in G$ so that $\langle g_j \rangle = \langle g_i \rangle^y$.  This implies that $\langle x_j N \rangle = \langle g_j N \rangle = \langle g_i N \rangle^{yN} = \langle x_i N \rangle^{yN}$.  Since the $x_i$'s were chosen to be representatives, we see that $i = j$.  Because the $\langle g_i \rangle$'s form a subset of a set of representatives of the conjugacy classes of maximal cyclic subgroups of $G$, we have $\eta (G/N) \le \eta (G)$. 

Now assume $\eta (G/N) = \eta (G)$.  I.e., we may assume that $\langle g_1 \rangle, \dots, \langle g_n \rangle$ are representatives of the maximal cyclic subgroups of $G$.  Suppose there exists $x \in N \setminus G^-$.  Then $\langle x \rangle$ is a maximal cyclic subgroup of $G$.  Notice that $g_1, \dots, g_n$ all lie outside of $N$.  Thus, $\langle x \rangle$ will not be conjugate to any $\langle g_i \rangle$, a contradiction.  Hence, we must have that $N \subseteq G^-$.  

We now work to show $(G/N)^- = \{gN \mid gN \subseteq G^-\}$.  Suppose that $xN \in (G/N)^-$ and suppose there exists $y \in xN \setminus G^-$.  It follows that $\langle y \rangle$ is maximal cyclic in $G$ and $\langle yN \rangle$ is not maximal cyclic in $G/N$.  In particular, $\langle y \rangle$ cannot be conjugate to one of the $\langle g_i \rangle$'s, contradiction.  Thus, we have $xN \subseteq G^-$.  Conversely, suppose now that $xN \subseteq G^-$.  We know that if $xN \not\in (G/N)^-$, then $\langle xN \rangle$ is maximal cyclic in $G/N$.  We have seen $\langle xN \rangle$ is conjugate to $\langle g_i N \rangle$ for some $i$, and so, there exists $g \in G \setminus G^-$ so that $xN = gN$, and we do not have $xN \subseteq G^-$, a contradiction.  Thus, we have that $(G/N)^- = \{ gN \mid gN \subseteq G^-\}$.  

We now prove that if $x \in G \setminus G^-$ and $y \in xN \setminus G^-$, then $y$ is conjugate to a generator of $\langle x \rangle$.  We have that $\langle x \rangle$ and $\langle y \rangle$ are maximal cyclic subgroups of $G$ and $\langle xN \rangle =  \langle yN \rangle$ is a maximal cyclic subgroup of $G/N$.  Using the fact that $\eta (G) = \eta (G/N)$, this implies that $\langle x \rangle$ and $\langle y \rangle$ must be conjugate to the same $\langle g_i \rangle$ and so, they are conjugate to each other.  We conclude that $y$ must be conjugate to a generator of $\langle x \rangle$.

Conversely, assume that $N \subseteq G^-$, $(G/N)^- = \{ gN \mid gN \subseteq G^- \}$ and every element $x \in G \setminus G^-$ satisfies the condition that every element of $xN \setminus G^-$ is conjugate to a generator of $\langle x \rangle$.  Define a map from the maximal cyclic subgroups of $G$ to the maximal cyclic subgroups of $G/N$ by $C \mapsto CN/N$.  Suppose $C$ is a maximal cyclic subgroup of $G$.  Then $C =  \langle x \rangle$ for some $x \in G \setminus G^-$.  It follows that $xN \in G/N \setminus (G/N)^-$.  We see that $CN/N =  \langle xN \rangle$ is a maximal subgroup of $G/N$.  Also, if $C = \langle x \rangle = \langle y\rangle$, then $\langle xN \rangle = \langle yN \rangle$, so the map is well-defined.  Suppose $D$ is a maximal cyclic subgroup of $G/N$.  Then $D =  \langle xN \rangle$ where $xN \in G/N \setminus (G/N)^-$.  As we have seen, this implies that $xN \setminus G^-$ contains some element $g$, and thus, $\langle g \rangle$ is a maximal cyclic subgroup of $G$.  Also, $\langle g \rangle N/N = \langle gN \rangle = \langle xN \rangle = D$.  This implies that the map is onto.  Finally, suppose that $x, y \in G \setminus G^-$ yields $\langle xN \rangle = \langle yN \rangle$.  This implies that some generator of $\langle y \rangle$  lies in $xN$.  By hypothesis, this generator will be conjugate to a generator of $\langle x \rangle$.  This implies that $\langle x \rangle$ and $\langle y \rangle$ are conjugate.  Hence, our map yields a bijection of conjugacy classes, and we conclude that $\eta (G) = \eta (G/N)$.

We now suppose that $G$ is a $p$-group and assume $N$ is a normal subgroup so that $\eta (G/N) = \eta (G)$.  Suppose $x \in G \setminus G^{\{p\}} = G \setminus G^-$.  The above work shows that $xN \in G/N \setminus (G/N)^{\{p\}}$.  Let $y \in xN$ and suppose $y \in G^{\{p\}}$, so $y = g^p$ for some element $g \in G$.  This implies that $xN = yN = g^pN = (gN)^p$ and this contradicts $xN \in G/N \setminus (G/N)^{\{p\}}$.  We conclude that $xN \cap G^-$ is empty.  This implies that $G \setminus G^-$ and $G^-$ are both unions of cosets of $G$.  

If $G^-$ is a union of $N$-cosets and $\eta (G/N) = \eta (G)$, then when $x \in G \setminus G^-$, we must have that $xN \cap G^-$ is empty, and so, every element $y \in xN$ lies in $G \setminus G^-$ and so applying conclusion (2)(c), we have that $y$ is conjugate to a generator of $\langle x \rangle$.  Conversely, suppose for every $x \in G \setminus G^-$ that every element of $xN$ is conjugate to a generator of $\langle x \rangle$.  Consider the element $n \in N \setminus \{ 1 \}$.  Since $\langle n \rangle$ cannot be conjugate to $\langle 1 \rangle$, we must have that $n \in G^-$, and so, $N \subseteq G^-$.  If $x \in G \setminus G^-$, then every generator of $\langle x \rangle$ will lie in $G \setminus G^-$.  Since every element of $xN$ is conjugate to a generator of $\langle x \rangle$, we have that $xN \subseteq G \setminus G^-$, and we deduce that $G\setminus G^-$ and hence, $G^-$ is a union of $N$-cosets.  It follows that $G^- N \subseteq G^-$, and we deduce that $G^- N = G^-$.  This implies that $(G/N)^- = \{ gN \mid gN \subseteq G^-\} = \{ gN \mid g \in G^- \} = G^-N/N$.  It follows that (a), (b), and (c) of (2) are met, so $\eta (G) = \eta (G/N)$.  Notice that we have proved both (4) and (5).
$\Box$\\[1ex]

We note that we cannot completely drop the hypothesis that $G$ is a $p$-group in conclusion (3).  Let $G$ be the dihedral group of order $30$.  Let $M$ be the normal subgroup of order $3$ and $N$ the normal subgroup of order $5$.  It follows that $G/M$ is a dihedral group of order $10$ and $G/N$ is a dihedral group of order $6$.  We see that $\eta (G/M) = \eta (G/N) = \eta (G) = 2$.  On the other hand, it is not difficult to see that $G^- = N \cup M$, and $G^-$ is not a union of $N$-cosets or $M$-cosets.  If $x \in NM \setminus G^-$, then it is not difficult to see that neither $xN \cap G^-$ nor $xM \cap G^-$ is empty.  Also, $G^- N = NM = G^-M$.

In Section \ref{prod sect}, we considered the case where $G$ is the semi-direct product of $H$ acting on $N$.  Using Lemma \ref{centre}, we saw that $\eta (G) \ge \eta^* (N)$.  Since $H \cong G/N$, we can use Theorem \ref{quot} (1) to see that $\eta (G) \ge \eta (G/N) = \eta (H)$.  Notice that the examples of $D_8$ acting on $Z_3 \times Z_3$ is an example where $\eta (G) = \eta (H)$.

For $p$-groups, we can find a largest normal subgroup such that $\eta$ of the quotient is equal to $\eta$ of the group.  In this situation, we also cannot completely drop the hypothesis that $G$ is a $p$-group from Corollary \ref{products}.  Again, let $G$ be the dihedral group of order $30$, let $M$ be the normal subgroup of order $3$, and let $N$ the normal subgroup of order $5$.  We have seen that $\eta (G/M) = \eta (G/N) = \eta (G) = 2$.  Notice that $G/NM$ is cyclic of order $2$, and so, $\eta (G/MN) = 1$.

\begin{corollary}\label{products}
Let $G$ be a noncyclic finite $p$-group for some prime $p$.
\begin{enumerate}
\item [(i)] Let $N$ and $M$ be normal subgroups of $G$.  If $\eta (G) = \eta(G/N) = \eta (G/M)$, then $\eta (G) = \eta (G/NM)$.
\item [(ii)] There exists a characteristic subgroup $X$ of $G$ so that $\eta (G) = \eta (G/X)$ and if $N$ is a normal subgroup $G$ so that $\eta (G/N) = \eta (G)$, then $N \le X$.
\end{enumerate}
\end{corollary}

{\bf Proof.}
(i) As $N$ and $M$ both lie inside $G^-$ by Theorem \ref{quot}, it follows that $NM \subseteq G^-$ again by Theorem \ref{quot}. Suppose now that $g \in G \setminus G^-$.  We need to show that every element in $gMN$ is conjugate to a generator of $\langle g \rangle$.  An arbitrary element of $gMN$ has the form $gmn$ where $m \in M$ and $n \in N$.  Applying Theorem \ref{quot} to $M$ in $G$, we see that every element of $gM$ is conjugate to a generator of $\langle g \rangle$.  In particular, $\langle gm \rangle$ is conjugate to $\langle g \rangle$.  This implies that $gm \in G \setminus G^-$.  We then apply Theorem \ref{quot} to $N$ in $G$ to see that every element in $gmN$ is conjugate to a generator of $\langle gm \rangle$.  We obtain $\langle gmn \rangle$ is conjugate to $\langle gm \rangle$.  Using the transitivity of conjugate subgroups, we see that $\langle g \rangle$ and $\langle gmn \rangle$ are conjugate.  We conclude that $gnm$ is conjugate to a generator of $\langle g \rangle$ as desired.   We have now shown that $NM$ satisfies the conditions of Theorem \ref{quot}, and so, $\eta (G/MN) = \eta (G)$.

(ii)  Let $\mathcal {M}  = \{ M \unlhd G \mid \eta (G/M) = \eta (G) \}$.  Let $X = \prod_{M \in \mathcal {M}} M$.  It is not difficult to see that $X$ is characteristic in $G$ and if $N$ is normal with $\eta (G/N) = \eta(G)$, then $N \le X$.  Since $|G|$ is finite, we can use (i) inductively to see that $\eta (G/X) = \eta (G)$.
$\Box$\\ 





Our next corollary states that $N \le G'$ when $\eta (G) = \eta (G/N)$ and every Sylow subgroup of $G/G'$ is noncyclic.  Note that if $G$ is a noncyclic $p$-group, then $G/G'$ is noncyclic.  It follows that this next corollary applies to $p$-groups.   

We now show that the hypothesis that $G/G'$ is not cyclic is necessary.  Let $M$ be a cyclic group of order $3$ and let $H$ be a cyclic group of order $4$.  Note that $H$ has a nontrivial action on $M$ and let $G$ be the resulting semidirect product. Take $N$ to be the subgroup of order $2$ in $H$ and note that $N = Z(G)$.  It is easy to see that $M = G'$.  Note that $G/N \cong S_3$, and so, $\eta (G/N) = 2$.  We claim that $\eta (G) = 2$, also.  To see this, observe that $MN$ is cyclic of order $6$.  Suppose $g \in G \setminus MN$.  It follows that $gN \in G/N \setminus MN/N$.  Since $G/N$ is isomorphic to $S_3$, we see that $gN$ is conjugate to $hN$ where $h$ is the generator of $H$.  Since $N \le H$, this implies that $g$ is conjugate to either $h$ or $h^3$.  We conclude that $\eta (G) =2$.  Finally, note that $N$ is not contained in $G' = M$.  We now present a second example to see that we need to assume that all Sylow subgroups of $G/G'$ are noncyclic, take $K$ to be a group so that $K/K'$ is noncyclic, and take $C$ to be a cyclic $p$-group for a prime $p$ so that $p$ does not divide $|K|$.  Let $G = K \times C$.  We see that $\eta (G/C) = \eta (K)= \eta (K \times C) = \eta (G)$ and $C$ is not contained in $G'$.  This is Theorem \ref{3rd main}.
 
\begin{corollary}\label{derived} 
Suppose $G$ is a group such that every nontrivial Sylow subgroup of $G/G'$ is not cyclic and $N$ a normal subgroup such that $\eta(G) = \eta(G/N)$. Then $N \leq G'$.
\end{corollary} 

{\bf Proof.}
Let $\Phi/G' = \Phi (G/G')$, so $\Phi$ is the preimage of the Frattini subgroup of $G/G'$.  We claim that $N \le \Phi$.  To obtain a contradiction, we suppose that $N$ is not contained in $\Phi$.  Let $M = N \cap \Phi$, and we have $M < N$.  We know that $G/\Phi$ is the direct product of elementary abelian $p$-groups for various primes $p$.  In particular, there is a subgroup $H$ so that $G = H(N\Phi)$ and $H \cap N\Phi = \Phi$.  This implies that $\Phi \le H$ and so, $G = N H$.  Also, $H \cap N = H \cap (N \Phi \cap N) = (H \cap N\Phi) \cap N = \Phi \cap N = M$.  It follows that $G/M = H/M \times N/M$.  Observe that $G/N \cong H/M$ and $G/H \cong N/M$.  We see that $N/M$ is nontrivial.  If $N/M$ is not cyclic, then by Lemma \ref{dir product} (i), we have $\eta (G/M) \ge (\eta (N/M))(\eta (H/M)) > \eta (H/M)$.   On the other hand, suppose $N/M$ is cyclic, and let $p$ be a prime that divides $|N:M|$.  Observe that $N/M \cong N\Phi/\Phi$ and $G/\Phi = N\Phi/ \Phi \times H/\Phi$.  Since we are assuming no Sylow subgroup of $G/G'$ is cyclic, we see that no Sylow subgroup of $G/\Phi$ is cyclic.  It follows that $p$ must divide $|H:\Phi|$, and so, $p$ divides $|H:M|$.  We now apply Lemma \ref{dir product} (iv) to see that $\eta (G/M) > \eta (H/M)$.  In both cases, we have proved that $\eta (G/M) > \eta (H/M)$.  By Theorem \ref{quot}, we have $\eta (G) \ge \eta (G/M) > \eta (H/M) = \eta (G/N)$.   This however contradicts $\eta (G) = \eta (G/N)$.  Thus, we must have that $N \le \Phi$.

Suppose $N$ is not contained in $G'$.  Then we can find $n \in N \setminus G'$ so that $nG'$ has order $p$ in $G/G'$ for some prime $p$.  Since $G/G'$ is not cyclic, we know that $G/\Phi$ is not cyclic.  Hence, we can find $g \in G \setminus \Phi$, so that $\langle nG' \rangle$ is not contained in $\langle gG' \rangle$.  In particular, $\langle gG' \rangle$ and $\langle gnG' \rangle$ are distinct subgroups of $G/G'$.   Notice that $g, gn \in G \setminus \Phi \subseteq G \setminus G^-$.  By Theorem \ref{quot}, we know that $\eta (G) = \eta (G/N)$ implies that every element in $gN \setminus G^-$ is conjugate to a generator of $\langle g \rangle$.  This implies that $gn$ is conjugate to a generator of  $\langle g \rangle$.  Hence, $gnG'$ is conjugate to an element of $\langle gG' \rangle$.  Since $G/G'$ is abelian and $\langle gnG' \rangle$ and $\langle gG' \rangle$ are distinct subgroups of $G/G'$, this is a contradiction.
$\Box$\\

Since it is well known when $G$ is a noncyclic $p$-group that $G/G'$ is not cyclic, Corollary \ref{derived} yields:

\begin{corollary}
If $G$ is a noncyclic $p$-group and $N$ is a normal subgroup so that $\eta (G/N) = \eta (G)$, then $N \le G'$.
\end{corollary}

Applying Theorem \ref{quot} to $p$-groups of exponent $p$, we see that $\eta (G)$ grows proportionally to $\log_p (|G|)$.

\begin{corollary}\label{exp}
Suppose $n \geq 2$ and $G$ is a $p$-group of order $p^n$ and exponent $p$.  Then $\eta(G) \geq n + p -1$.
\end{corollary}

{\bf Proof.} 
If $|G| = p^2$ then $\eta(G) = p+1$. We proceed by induction on the order of $G$ and assume $|G| \geq p^3$. Suppose $z$ is a central element of $G$, then $\eta (G/\langle z \rangle ) < \eta(G)$ by Theorem \ref{quot}. Furthermore, $G/\langle z \rangle$ is a noncyclic $p$-group of exponent $p$. The result follows.$\Box$\\

A covering $\{ H_i \}$ of a group $G$ is called a {\it partition} if $H_i \cap H_j$ is trivial for all $i \neq j$.   (I.e., $H_i \cap H_j = 1$ for $i \ne j$.)  Also, a group $G$ is {\it tidy} if for all $x \in G$ the set ${\rm Cyc}_G(x) = \{ y \in G: \langle x,y \rangle \;\;{\rm is} \;{\rm cyclic}\}$ is a subgroup of $G$.  Suppose $G$ is a finite $p$-group and $G$ has a covering that consists of cyclic groups that is also a partition.  Clearly, a cyclic group is tidy, so it follows from \cite[Corollary 2.5]{erfanian} that $G$ is tidy. The tidy $p$-groups have been classified in \cite{lewis}; they are either (i) cyclic, (ii) have exponent $p$ or (iii) $p=2$ and $G$ is dihedral or generalized quaternion.  Hence, if $G$ is a partitioned $p$-group, then either $G$ is cyclic and $\eta (G) = 1$, $G$ is dihedral or generalized quaternion and $\eta (G) = 3$, or $G$ has exponent $p$ and $\eta (G) \ge \log_p (|G|) + p - 1$.

The following result gives a dichotomy of normal subgroups.  Our proof requires that $G$ be a $p$-group, but we do not have any counterexamples when $G$ is not a $p$-group, so it may be possible to weaken the hypothesis on this corollary.  Finally, we have Theorem \ref{4th main}.

\begin{corollary} \label{eitheror}
Let $G$ be a noncyclic $p$-group and $N$ a nontrivial normal subgroup of $G$ such that $\eta(G) = \eta(G/N)$.
\begin{enumerate}
\item [(i)] Suppose $M \unlhd G$ then either $N \leq M$ or $M \leq G^-$.
\item [(ii)] If there exists a maximal cyclic group $\langle x \rangle$ which is normal, then $N \leq \langle x \rangle$; so $N$ is cyclic.
\end{enumerate}
\end{corollary}

{\bf Proof.} 
(i) Suppose that $N$ is not contained in $M$ and $M$ is not contained in
$G^-$.  Then we can find $n \in N \setminus M$ and $m \in M \setminus
G^-$.  It follows that $m \in M$ and $mn$ is not in $M$.  Since $M$ is normal,
$mn$ cannot be conjugate to any generator of $\langle m \rangle$.   Since $m$ is not in
$G^-$, we should have that every element of $mN$ is conjugate to a
generator of $\langle m \rangle$, and $mn$ is not; so we have a contradiction to Theorem \ref{quot}.

(ii) We are assuming that $\langle x \rangle$ is normal in $G$ and $x \in G \setminus G^-$.  Thus, $\langle x \rangle$ is not contained in $G^-$.  So by (i) we have $N \le \langle x \rangle$.
$\Box$\\[1ex]

\noindent Mariagrazia Bianchi:
 Dipartimento di Matematica F. Enriques,
 Universit\`a degli Studi di Milano, via Saldini 50,
20133 Milano, Italy.\\
mariagrazia.bianchi@unimi.it\\[1ex]
Rachel D. Camina: Fitzwilliam College, Cambridge, CB3 0DG, UK.\\
rdc26@cam.ac.uk\\[1ex]
Mark L. Lewis:  Department of Mathematical Sciences, Kent State University, Kent, Ohio, 44242 USA.\\
lewis@math.kent.edu\\[1ex]
Emanuele Pacifici: Dipartimento di Matematica e Informatica "U. Dini" (DIMAI),
Universit\'a degli Studi di Firenze,
Viale Morgagni 67/A, 50134 Firenze, Italy.\\
emanuele.pacifici@unifi.it\\[1ex]

\end{document}